\documentclass{article}
\usepackage{amsmath}
\usepackage{amsthm}
\usepackage{amssymb}
\usepackage{bm}
\usepackage[dvips]{graphicx}
\usepackage{latexsym}
\usepackage{mathtools}
\usepackage{geometry}
 \renewcommand{\equation}

\newtheorem{prop}{Proposition}[section]

\newtheorem{thm}{Theorem}[section]

\newtheorem{fact}{Fact}[section]

\theoremstyle{definition}
\newtheorem{rmk}{Remark}[section]
\newtheorem{defini}{Definition}[section]

\title{The subgroup of the Goeritz group of the Heegaard splitting induced by an openbook decomposition consisting of elements preserving the binding}
\author{Nozomu Sekino}
\date{}
\begin{document}
\maketitle

\begin{abstract}
When a 3-manifold admits an openbook decomposition, we get a Heegaard splitting by thickening a page. 
This splitting surface has a special multi curves coming from the binding. 
In this paper, we consider the subgroup of the Goeritz group of this Heegaard splitting, which is the mapping class group of the 3-manifold preserving the given Heegaard splitting, consisting of elements preserving the binding. 
This subgroup turned out to be the quotient of the subgroup of the orientation preserving mapping class group consisting of elements commuting with the monodromy by the subgroup generated by the Dehn twists along the boundary curves. 
We also get a criterion for the existence of an element of the Goeritz group which fixes the binding as a set and reverses the orientation. 
At last, we give some example of computation of a Goeritz group. 
\end{abstract}

\section{Introduction} \label{intro}
It is known that every connected orientable closed 3-manifold $X$ admits a {\it Heegaard splitting} $X=V^{+}\cup_{\Sigma}V^{-}$: 
This represents a connected orientable closed surface $\Sigma$ in $X$ and two handlebodies $V^{+}$ and $V^{-}$ bounded by $\Sigma$ in $X$ ``inner'' and ``outer'', respectively. 
The {\it Hempel distance} \cite{hempel} of a Heegaard splitting $X=V^{+}\cup_{\Sigma}V^{-}$ is defined to be the distance between the sets of the boundaries of properly embedded essential disks of $V^{+}$ and $V^{-}$ in (the 1-skelton of) the curve complex of $\Sigma$, and this measures some complexity of the Heegaard splitting. 

For a Heegaard splitting $X=V^{+}\cup_{\Sigma}V^{-}$, the set of isotopy classes of orientation preserving self-homeomorphisms of $X$ fixing $V^{+}$ (so also $V^{-}$ and $\Sigma$), where $\Sigma$ is fixed as a set during isotopies, forms a group by compositions. 
This group is called the {\it Goeritz group} of $X=V^{+}\cup_{\Sigma}V^{-}$. 
In other words, the Goeritz group of $X=V^{+}\cup_{\Sigma}V^{-}$ is a subgroup of the mapping class group of $\Sigma$ consisting of the elements whose representatives can extend to self-homeomorphisms of both of $V^{+}$ and $V^{-}$. 

About Goeritz groups, Johnson \cite{johnson1} showed that if the Hempel distance of a Heegaard splitting is at least $4$, then the Goeritz group of the Heegaard splitting is finite. 
This means that a high Hempel distance Heegaard splitting is very ``rigid''. 
On the other hand, the Goeritz group of a Heegaard splitting has an element of infinite order if the Hempel distance of the Heegaard splitting is at most $1$ (, twist along a reducing sphere or so-called ``eye glass twist'' using weakly reducible pair). 
In this case, determining the structure of the Goeritz group is difficult in general. 
It is unknown whether the Goeritz group of the Heegaard spliting of genus greater than $3$ of the $3$-sphere is finitely generated or not. 
By the sequence of works \cite{goeritz},\cite{scharlemann},\cite{cho1},\cite{cho2},\cite{ck1},\cite{ck2},\cite{ck3},\cite{ck4}, 
a finite presentation of the Goeritz group of every Heegaard splitting of genus two of Hempel distance at most $1$ is known. 
Freedman and Scharlemann \cite{fs} gave a finite generating set of the Goeritz group of the Heegaard splitting of genus three of $3$-sphere. 
For the case where Hempel distance is $2$ or $3$, the Goeritz groups may be finite or infinite. 
For examples of infinite cases, the Goeritz group of a Heegaard splitting induced by an open book decomposition of some type  \cite{johnson2} or by a twisted book decomposition of some type \cite{ik1} is infinite, these examples have Hempel distance $2$. 
For examples of finite cases, the Goeritz group of a Heegaard splitting having some ``keen''-type property, which represents the rigidity of a Heegaard splitting \cite{ik1} is finite, and a concrete example are given in \cite{sekino}, these examples have Hempel distance at least $2$. 
About finitely generativity, the Goeritz group of a Heegaard splitting having ``thick isotopy'' property and finiteness of the mapping class group of the ambient manifold \cite{iguchi} is finitely generated, these examples have Hempel distance at least $2$. 
The example of a infinite Goeritz group of a Heegaard splitting of Hempel distance $3$ is not found for now. 

In this paper, we focus on the Goeritz groups of the Heegaard splittings induced by openbook decompositions. 
As stated before, their Hempel distance are less than or equal to $2$, and possibly gives an element of infinite order. 
An openbook decomposition of a 3-manifold is a way of presentation of the manifold as some quotient of oriented compact surface bundle over a circle. The image of a surface under the quotient is called a {\it page}. 
It is known that every connected oriented closed 3-manifold admits an openbook decomposition \cite{alexander}. 
By thickening a page, we get a Heegaard splitting, whose splitting surface has some special multi curve, called a {\it binding}, coming from the boundary of a page. 
There is an isotopy of the manifold whose endpoints homeomorphisms preserve the Heegaard decomposition: Sliding handlebodies one rap along the circle of the fiber structure. 
This map fixes the binding. 
We consider the subgroup of the Goeritz group consisting of the elements preserving the binding. 
The following is a result. Notations, conventions and terminologies are in the last of Section~\ref{intro} and Section~\ref{terminologies}.

\begin{thm} \label{thm}
Let $M$ be a connected closed oriented 3-manifold with an openbook decomposition each of whose page is homeomorphic to $\Sigma_{g,b}$ and whose monodromy is $\phi$, and $M=V\cup W$ the Heegaard splitting induced by the openbook decomposition. 
Take any orientation reversing involution $\iota$ of $\Sigma_{g,b}$. 
Then
\begin{itemize}
\item $\mathcal{G}_{bind}\cong \{ [f]\in MCG_{+}(\Sigma_{g,b})\mid [\phi \circ f]= [f\circ \phi] \}/ \left< [t_{c}]\mid  c {\rm \ is\ a\ boundary\ component\ of\ }\Sigma_{g,b} \right>$. 
\item $\mathcal{G}$ has an element one of whose representatives fixes the binding as a set and reverses the orientation of the binding if and only if there exists $[f]\in MCG_{+}(\Sigma_{g,b})$ such that $\phi^{-1}\circ f=f\circ \iota \circ \phi \circ \iota$.
\end{itemize}
\end{thm}

The rest of the paper is organized as follows. 
In Section~\ref{terminologies}, we recall terminologies used in the title. We give a Heegaard diagram we will work on. 
In Section~\ref{tools}, we list some operations needed for proofs. 
In Section~\ref{formerproof}, we prove the former part of Theorem~\ref{thm}. For a given map on the Heegaard surface, we search the condition for its extending to each of the handlebodies. 
In Section~\ref{latterproof}, we prove the latter part of Theorem~\ref{thm}. We assume that there exists a binding revering map, and search what it should satisfy. 
In Section~\ref{example}, we give some example of a Goeritz group of the Heegaard splitting induced by an openbook decomposition whose monodromy is sufficiently complex so that every element of the Goeritz group fixes the binding as a set. 

\subsection*{Notations and convenions}
In this paper, $\Sigma_{g,b}$ denotes the connected oriented compact surface of genus $g$ with $b$ boundary components. 
In addition, we suppose $b\geq1$ throughout the paper. 
For an oriented connected 1-manifold $I$ with boundary, the starting point and the terminal point of $I$ are denoted by $s(I)$ and $t(I)$, respectively. 
The 1-manifold obtained by reversing the orientation of $I$ is denoted by $rev(I)$. 
Note that $rev\left(rev(I)\right)=I$. 
For two real number $a$,$b$ with $a<b$, the interval $[a,b]$ is oriented so that $t([a,b])=\{b\}$. 
Words in some alphabets are read from the left to the right. 
For two word $\alpha$ and $\beta$, their commutator $\alpha \beta {\alpha}^{-1} {\beta}^{-1}$ is denoted by $[[\alpha, \beta]]$. 
The orientation preserving, boundary fixing mapping class group of $\Sigma_{g,b}$ is denoted by $MCG_{+}(\Sigma_{g,b})$. 
For a self-homeomorphism $f$ of a topological space, $[f]$ denotes the mapping class of $f$. 
The composition of maps $f$ and $g$ are denoted by $f\circ g$ and this is read from the right to the left.  
For an (oriented) simple closed curve $c$ in $\Sigma_{g,b}$, $t_c$ denotes the right-hand Dehn twist along $c$.

\subsection*{Acknowledgements}
The author would like to thank professor Sangbum Cho and professor Yuya Koda for introducing him the area of Goeritz groups.
%He also is grateful to referee for his or her kindness, patience and correcting many mistakes in the proofs and the arguments.

\section{The terminologies in the title} \label{terminologies}

When a connected oriented closed 3-manifold $M$ admits an openbook decomposition, 
$M$ is represented as the quotient of $\Sigma_{g,b}\times[0,2]$ for $b\geq1$ by the relations identifying $(p,2)$ with $\left( \phi(p),0 \right)$ for every $p\in \Sigma_{g,b}$ and $(q,t)$ with $(q,0)$ for every $q\in \partial \Sigma_{g,b}$ and $t\in[0,2]$, where $\phi$ is an orientation preserving homeomorphism fixing the boundary pointwise, called the {\it monodromy}. 
The image of $\partial \Sigma_{g,b}\times \{0\}$ is called the {\it binding}. 

Let $V$ and $W$ be the quotients of $\Sigma_{g,b}\times[0,1]$ and $\Sigma_{g,b}\times[1,2]$, respectively. 
Note that $V$ and $W$ are handlebodies of genus $2g+b-1$ and $M$ is obtained by pasting $V$ and $W$ along their boundaries. 
Thus we have a Heegaard splitting $M=V\cup W$. 
This splitting is called the {\it Heegaard splitting induced by the openbook decomposition}. 
Note that the binding of the openbook decomposition is on the splitting surface. 
%This multicurve on the splitting surface is also called the binding. 

\subsection{A Heegaard diagram of $M=V\cup W$} \label{hdiagram}
Let $\mathbb{J}=(J_{1},\dots,J_{2g+b-1})$ be any ordered tuple of pairwise disjoint properly embedded oriented arcs in $\Sigma_{g,b}$ such that they cut $\Sigma_{g,b}$ into a disk. 
Set 
\begin{itemize}
\item $S=\left( \Sigma_{g,b}\times \{0\}\right)\cup \left( \partial \Sigma_{g,b}\times[0,1]\right) \cup \left( \Sigma_{g,b} \times\{1\}\right)$, whose orientation coming from $\Sigma_{g,b}\times \{1\}$,  
\item $\mathbb{A}=(A_1,\dots,A_{2g+b-1})$, where $A_{i}=\left( J_{i}\times\{0\}\right) \cup \left( t(J_{i})\times [0,1]\right) \cup \left( rev(J_{i})\times \{1\}\right) \cup \left( s(J_{i})\times rev([0,1])\right)$, and
\item $\mathbb{B}=(B_1,\dots,B_{2g+b-1})$, where $B_{i}=\left( \phi(J_{i})\times\{0\}\right) \cup \left( t(J_{i})\times [0,1]\right) \cup \left( rev(J_{i})\times \{1\}\right) \cup \left( s(J_{i})\times rev([0,1])\right)$.
\end{itemize}
Then the triplet $(S; \mathbb{A},\mathbb{B})$ defines an (ordered) Heegaard diagram of $M=V\cup W$. 
Note that curves in both of  $\mathbb{A}$ and $\mathbb{B}$ are oriented, and cut $S$ into a connected planar surface. 
We regard $S$ as in $\mathbb{R}^{3}$. 
Note that the orientation of $\Sigma_{g,b}\times \{0\} \subset S$ coming from $S$ is the opposite of that coming from the orientation of $\Sigma_{g,b}$. 
Note also that $V$ and $W$ correspond to $S_{\mathbb{A}}$ and $S_{\mathbb{B}}$, respectively, where $S_{\mathbb{A}}$ and $S_{\mathbb{B}}$ denote the handlebodies obtained by pasting $2$-handles along curves in $\mathbb{A}$ inside of $S$ and $\mathbb{B}$ outside of $S$, respectively, and filling the resulting spheres with $3$-balls. 
Note also that the orientation of $S$ is same as that coming from handlebody $S_{\mathbb{A}}$ and is opposite to that coming from handlebody $S_{\mathbb{B}}$, where handlebodies are oriented as subsets in $\mathbb{R}^3$. 
In this diagram, the binding corresponds to $\partial \Sigma_{g,b}\times \{\frac{1}{2}\}\subset S$. 

\begin{rmk}
When $2g+b-1\geq2$, there is a simple closed curve $l$ in $\Sigma_{g,b}$ which is disjoint from $J_1$. 
Then $l\times \{1\}\subset S$ is disjoint from $A_1$ and $B_1$ and we see that the Hempel distance of the Heegaard splitting is less than or equal to $2$. 
\end{rmk}

\subsection{The subgroup of the Goeritz group of $M=V\cup W$ consisting of elements preserving the binding} \label{subgroup}

The Goeritz group $\mathcal{G}$ of $M=V\cup W$ is the group consisting of isotopy classes of orientation preserving self-homeomorphism of $M$ preserving $V$ as a set (and also $W$) where isotopies preserve $\partial V$ as a set. 
Using a diagram $(S;\mathbb{A},\mathbb{B})$, we regard $\mathcal{G}$ as the subgroup of $MCG_{+}(S)$ consisting of elements preserving $S_{\mathbb{A}}$ and $S_{\mathbb{B}}$ i.e. $[h]\in MCG_{+}(S)$ such that $h(A_{i})$ and $h(B_{i})$ bounds disks in $S_{\mathbb{A}}$ and $S_{\mathbb{B}}$, respectively for all $i=1,\dots,2g+b-1$. 
We consider a subgroup of $\mathcal{G}$ consisting of elements preserving the binding. 
\begin{defini} \label{bindpres}
Let $\mathcal{G}_{bind}$ denote the subgroup of $\mathcal{G}$ consisting of elements preserving the binding i.e. one of each representatives fixes the binding pointwise. 
\end{defini}

\begin{rmk}
Of course, not all elements in $\mathcal{G}$ preserve the binding in general. 
However, when the monodromy is sufficiently ``complex'' so that the Heegaard splitting admits unique binding structure, $\mathcal{G}_{bind}$ is a normal subgroup of $\mathcal{G}$ of index less than or equal to $2$. 
The index is $2$ when there exists an element of $\mathcal{G}$, one of whose representatives fixes the binding as a set and reverses the orientation. 
Note that if the orientation of one component of the binding is reversed by an element of $\mathcal{G}$ fixing the binding as a set, then those of the other components are also reversed since the binding separates $S$ into two components. 
\end{rmk}

\section{Preparations for proofs}\label{tools}
In this section, we list definitions and remarks needed for proofs.
\subsection{The words of oriented arc on $\Sigma_{g,b}$} \label{arcword}
Let $\mathbb{J}=(J_{1},\dots,J_{2g+b-1})$ be an ordered tuple of pairwise disjoint properly embedded oriented arcs in $\Sigma_{g,b}$ such that they cut $\Sigma_{g,b}$ into a disk. 
For a properly embedded oriented arc $I$ in $\Sigma_{g,b}$ such that each of the endpoints of $I$ is not on arcs in $\mathbb{J}$, we assign the word $W_{\mathbb{J}}(I)$ in alphabets $\{x_{1},\dots,x_{2g+b-1}$ as follows: 
Follow $I$ from $s(I)$, and we start with the empty word. If we hit $I_{i}$ from the right side (or left side), then we add $x_{i}$ (or ${x_{i}}^{-1}$, respectively) to the right of the word we have. 

\begin{rmk}\label{uniqueword}
Let $\mathbb{J}$ be as the above. 
Let $I$ and $I'$ be two properly embedded oriented arcs on $\Sigma_{g,b}$ such that each of the endpoints is not on the arcs in $\mathbb{J}$ and that $t(I)=t(I')$ and $s(I)=s(I')$. 
Suppose that $W_{\mathbb{J}}(I)=W_{\mathbb{J}}(I')$. 
Then $I$ and $I'$ are isotopic in $\Sigma_{g,b}$ relative to the endpoints. 
This is because there is the unique way (up to isotopy) to realize the given word with fixed endpoints since the arcs in $\mathbb{J}$ cut $\Sigma_{g,b}$ into a disk. 
\end{rmk} 

\subsection{The cyclic words of oriented simple closed curves on the boundary of a handlebody} \label{sccword}
Let $H$ be a handlebody of genus $n\in \mathbb{N}$. 
Take an ordered tuple of pairwise disjoint properly embedded disks $\mathbb{D}=(D_1,\dots,D_{n})$ such that they cut $H$ into a ball. 
The ordered tuple of the oriented boundaries of disks in $\mathbb{D}$ is denoted by $\partial \mathbb{D}$. 
For an oriented simple closed curve $l$ in $\partial H$, we assign the cyclic word $\mathcal{W}_{\partial \mathbb{D}}(l)$ in alphabets $\{x_{1},\dots,x_{n}\}$ as follows: 
Follow $l$ from any point, and we start with the empty word. If we hit $\partial D_{i}$ from the right side (or left side), then we add $x_{i}$ (or ${x_{i}}^{-1}$, respectively) to the right of the word we have. 
When we return to the starting point, connect the last of the word to the initial of it to get a cyclic word. 
We will write cyclic words as if they were not cyclic because of space limitation. 
Note that $\mathcal{W}_{\mathbb{\partial \mathbb{D}}}(l)$ is the empty word if and only if $l$ bounds a disk in $H$. 

\subsection{Operations on the (cyclic) words of arcs and simple closed curves} \label{operations}
We list the operations on the (cyclic) words of arcs and simple closed curves. 

\begin{defini} \label{closure} (Closure)\\
For a word $W$ in some alphabets, $Cl(W)$ denotes the cyclic word obtained by connecting the last of $W$ to the initial of $W$. 
\end{defini}

\begin{defini} \label{reflection} (Reflection)\\
For a word $W$ in alphabet $\{x_{1},\dots,x_{n}\}$, $\overline{W}$ denotes the word obtained from $W$ by changing ${x_{i}}^{\pm1}$ to ${x_{i}}^{\mp1}$ in $W$ for all $i=1,\dots,n$. 
\end{defini}

%\begin{defini} \label{reversal} (Reversal)\\
%For an oriented arc $I$, $rev(I)$ denotes the oriented arc $I$ with the opposite orientation. 
%\end{defini}

\begin{rmk}\label{ref,rev,homeo}
Let $\mathbb{J}=(J_{1},\dots,J_{2g+b-1})$ be an ordered tuple of properly embedded pairwise disjoint oriented arcs in $\Sigma_{g,b}$ which cut $\Sigma_{g,b}$ into a disk. 
Let $I$ be a properly embedded oriented arc in $\Sigma_{g,b}$ such that whose endpoints are not on the arcs in $\mathbb{J}$. 
Let $f$ be a self-homeomorphism of $\Sigma_{g,b}$ fixing the boundary pointwise. Then

\begin{itemize}
\item $W_{f(\mathbb{J})}\left(f(I)\right)=W_{\mathbb{J}}(I)$. 
\item $W_{\mathbb{J}}\left( rev(I) \right)=\left( W_{\mathbb{J}}(I) \right)^{-1}$. 
\item $W_{rev(\mathbb{J})}(I)=\overline{W_{\mathbb{J}}(I)}$. 
\end{itemize}
, where $f(\mathbb{J})$ and $rev(\mathbb{J})$ denote $\left( f(J_1),\dots, f(J_{2g+b-1})\right)$ and $\left( rev(J_1),\dots, rev(J_{2g+b-1})\right)$, respectively.
\end{rmk}

\section{A proof of the former part of Theorem~\ref{thm}}\label{formerproof}
In this section, we give a proof of the former part of Theorem~\ref{thm}. 
Take  an ordered tuple of pairwise disjoint properly embedded oriented arcs in $\Sigma_{g,b}$ such that they cut $\Sigma_{g,b}$ into a disk, $\mathbb{J}=(J_{1},\dots,J_{2g+b-1})$. 

Take $[h]\in MCG_{+}(S)$ such that $h$ preserves the binding $\left( \partial \Sigma_{g,b}\times \{\frac{1}{2}\} \right) \subset S$ pointwise. 
We assume that $h$ fixes $\left( \partial \Sigma_{g,b}\times [0,1] \right) \subset S$ pointwise. 
In this situation, we get  two orientation preserving self-homeomorphisms of $\Sigma_{g,b}$ fixing the boundary pointwise $f_{00}$ and $f_{11}$ such that
\begin{itemize}
\item $h\mid_{\Sigma_{g,b}\times \{0\}}(p,0)=\left( f_{00}(p),0\right)$ for all $p\in \Sigma_{g,b}$.
\item $h\mid_{\Sigma_{g,b}\times \{1\}}(p,1)=\left( f_{11}(p),1\right)$ for all $p\in \Sigma_{g,b}$.
\end{itemize}

Let ${J_{i}}'$ be an oriented arc on $\Sigma_{g,b}$ obtained by slightly sliding $J_{i}$ to the right side so that the endpoints of ${J_{i}}'$ are not on the arcs in $\mathbb{J}$. 
And let ${A_{i}}'$ and ${B_{i}}'$ be oriented simple closed curves on $S$ obtained by replacing all $J_{i}$'s with ${J_{i}}'$'s in the definitions of $A_{i}$ and $B_{i}$, respectively. 
Note that ${A_{i}}'$ and ${B_{i}}'$ are parallel to $A_i$ and $B_i$, respectively. 

We consider the conditions for $h$'s preserving $S_{\mathbb{A}}$ in Subsection~\ref{former_sa} and $S_{\mathbb{B}}$ in Subsection~\ref{former_sb}. 
Then we define a surjective homomorphism from some subgroup of $MCG_{+}(\Sigma_{g,b})$ to $\mathcal{G}_{bind}$, and determine the kernel of the homomorphism in Subsection~\ref{homomorphism}. 

\subsection{A condition for preserving $S_{\mathbb{A}}$} \label{former_sa}
For all $i=1,\dots,2g+b-1$,\\
 $h({A_{i}}')=\left( f_{00}({J_{i}}')\times \{0\} \right)\cup \left( \{t({J_{i}}')\} \times[0,1] \right)\cup \left(  rev\left( f_{11}({J_{i}}')\right) \times \{1\} \right)\cup \left( \{s({J_{i}}')\} \times rev([0,1]) \right)$. \\
Then $\mathcal{W}_{\mathbb{A}}\left( h({A_{i}}') \right)=Cl\left( \overline{W_{\mathbb{J}}\left( f_{00}({J_{i}}') \right)} \cdot W_{rev(\mathbb{J})}\left( rev\left( f_{11}({J_{i}'}) \right)\right) \right)$. 
Note that the right side of $J_{i}\times \{0\}$ in $\Sigma_{g,b}\times\{0\}$ with the orientation coming from that of $\Sigma_{g,b}$ corresponds to the left side of $A_i$ in $\partial S_{\mathbb{A}}$, and that the right side of $rev(J_{i})\times \{1\}$ in $\Sigma_{g,b}\times\{1\}$ with the orientation coming from that of $\Sigma_{g,b}$, which is the same as that coming from that of $S$, corresponds to the right side of $A_i$ in $\partial S_{\mathbb{A}}$. 
Using Remark~\ref{ref,rev,homeo}, we know that $\overline{W_{\mathbb{J}}\left( f_{00}({J_{i}}') \right)}=W_{rev(\mathbb{J})}\left( f_{00}({J_{i}}') \right)$ and that $W_{rev(\mathbb{J})}\left( rev\left( f_{11}({J_{i}'}) \right)\right)=\left( W_{rev(\mathbb{J})}\left( f_{11}({J_{i}}') \right) \right)^{-1}$. 
Therefore, $h({A_{i}}')$ bounds a disk in $S_{\mathbb{A}}$ (i.e $\mathcal{W}_{\mathbb{A}}\left( h({A_{i}}') \right)$ is the empty word) if and only if $W_{rev(\mathbb{J})}\left( f_{00}({J_{i}}') \right)=W_{rev(\mathbb{J})}\left( f_{11}({J_{i}}') \right)$. 
By Remark~\ref{uniqueword} since $f_{00}$ and $f_{11}$ fix the endpoints of ${J_{i}}'$, this is equivalent to that $f_{00}({J_{i}'})$ and $f_{11}({J_{i}}')$ are isotopic in $\Sigma_{g,b}$ relative to their endpoints. 

Since curves in $\mathbb{A}$ cut $S$ into a planar surface, $h$ preserves $S_{\mathbb{A}}$ if and only if $h({A_{i}}')$ bounds a disk in $S_{\mathbb{A}}$ for all $i=1,\dots,2g+b-1$. 
By the argument in previous paragraph, it is equivalent to that $f_{00}({J_{i}'})$ and $f_{11}({J_{i}}')$ are isotopic in $\Sigma_{g,b}$ relative to their endpoints for all $i=1,\dots,2g+b-1$. 
Note that we can isotope $f_{00}({J_{i}'})$ to $f_{11}({J_{i}}')$ for all $i$ simultaneously. 
Since arcs in $\mathbb{J}$ cut $\Sigma_{g,b}$ into a disk, this is equivalent to that $f_{00}$ and $f_{11}$ are isotopic under the isotopy fixing the boundary.

\subsection{A condition for preserving $S_{\mathbb{B}}$} \label{former_sb}
For all $i=1,\dots,2g+b-1$,\\
 $h({B_{i}}')=\left( f_{00}\left( \phi({J_{i}}')\right)\times \{0\} \right)\cup \left( \{t({J_{i}}')\}\times[0,1] \right)\cup \left(  rev\left( f_{11}({J_{i}}')\right) \times \{1\} \right)\cup \left( \{s({J_{i}}')\}\times rev([0,1]) \right)$. \\
Then $\mathcal{W}_{\mathbb{B}}\left( h({B_{i}}') \right)=Cl\left( W_{\phi(\mathbb{J})}\left( f_{00}\left( \phi({J_{i}}')\right) \right) \cdot \overline{W_{rev(\mathbb{J})}\left( rev\left( f_{11}({J_{i}'}) \right)\right) }\right)$. 
Note that the right side of $\phi(J_{i})\times \{0\}$ in $\Sigma_{g,b}\times\{0\}$ with the orientation coming from that of $\Sigma_{g,b}$ corresponds to the right side of $B_i$ in $\partial S_{\mathbb{B}}$, and that the right side of $rev(J_{i})\times \{1\}$ in $\Sigma_{g,b}\times\{1\}$ with the orientation coming from that of $\Sigma_{g,b}$, which is the same as that coming from that of $S$, corresponds to the left side of $B_i$ in $\partial S_{\mathbb{B}}$. 
Using Remark~\ref{ref,rev,homeo}, we know that $W_{\phi(\mathbb{J})}\left( f_{00}\left( \phi({J_{i}}')\right) \right)=W_{\mathbb{J}}\left( \phi^{-1}\circ f_{00} \circ \phi ({J_{i}}') \right) $ and that $\overline{W_{rev(\mathbb{J})}\left( rev\left( f_{11}({J_{i}'}) \right)\right)}=\left( W_{\mathbb{J}}\left( f_{11}({J_{i}}') \right) \right)^{-1}$. 
Therefore, $h({B_{i}}')$ bounds a disk in $S_{\mathbb{B}}$ (i.e $\mathcal{W}_{\mathbb{B}}\left( h({B_{i}}') \right)$ is the empty word) if and only if $W_{\mathbb{J}}\left( \phi^{-1}\circ f_{00}\circ \phi({J_{i}}') \right)=W_{\mathbb{J}}\left( f_{11}({J_{i}}') \right)$. 
By Remark~\ref{uniqueword} since $f_{00}$, $f_{11}$ and $\phi$ fix the endpoints of ${J_{i}}'$, this is equivalent to that $\phi^{-1}\circ f_{00}\circ \phi({J_{i}'})$ and $f_{11}({J_{i}}')$ are isotopic in $\Sigma_{g,b}$ relative to their endpoints. 

Since curves in $\mathbb{B}$ cut $S$ into a planar surface, $h$ preserves $S_{\mathbb{B}}$ if and only if $h({B_{i}}')$ bounds a disk in $S_{\mathbb{B}}$ for all $i=1,\dots,2g+b-1$. 
By the argument in previous paragraph, it is equivalent to that $\phi^{-1}\circ f_{00}\circ \phi({J_{i}'})$ and $f_{11}({J_{i}}')$ are isotopic in $\Sigma_{g,b}$ relative to their endpoints for all $i=1,\dots,2g+b-1$. 
Note that we can isotope $\phi^{-1}\circ f_{00}\circ \phi({J_{i}'})$ to $f_{11}({J_{i}}')$ for all $i$ simultaneously. 
Since arcs in $\mathbb{J}$ cut $\Sigma_{g,b}$ into a disk, this is equivalent to that $\phi^{-1}\circ f_{00}\circ \phi$ and $f_{11}$ are isotopic under the isotopy fixing the boundary. 

\subsection{A surjective map to $\mathcal{G}_{bind}$ and its kernel} \label{homomorphism}
Let $Comm([\phi])$ denote the subgroup of $MCG_{+}(\Sigma_{g,b})$ consisting of elements commuting with $[\phi]$, $\{ [f]\in MCG_{+}(\Sigma_{g,b})\mid [\phi \circ f]= [f\circ \phi] \}$. 
\begin{defini}
For $f$ which is a representative of an element in $Comm([\phi])$, we define an orientation preserving self-homeomorphism of $S$, denoted by $F(f)$ as 
\begin{itemize}
\item $F(f)\mid_{\partial \Sigma_{g,b}\times[0,1]}(q,t)=(q,t)$ for every $q\in \partial \Sigma_{g,b}$ and $t\in[0,1]$, 
\item $F(f)\mid_{\Sigma_{g,b}\times \{0\}}(p,0)=\left(f(p),0\right)$ for every $p\in\Sigma_{g,b}$, and 
\item $F(f)\mid_{\Sigma_{g,b}\times \{1\}}(p,1)=\left(f(p),1\right)$ for every $p\in\Sigma_{g,b}$.
\end{itemize}
By the argument in previous two subsection, $F(f)$ preserves the binding, $S_{\mathbb{A}}$ and $S_{\mathbb{B}}$. 
This map induces a homomorphism $[F]: Comm([\phi])\longrightarrow \mathcal{G}_{bind}$ and we know that this is surjective by the argument in previous two subsection. 
\end{defini}

\begin{figure}[htbp]
 \begin{center}
  \includegraphics[width=100mm]{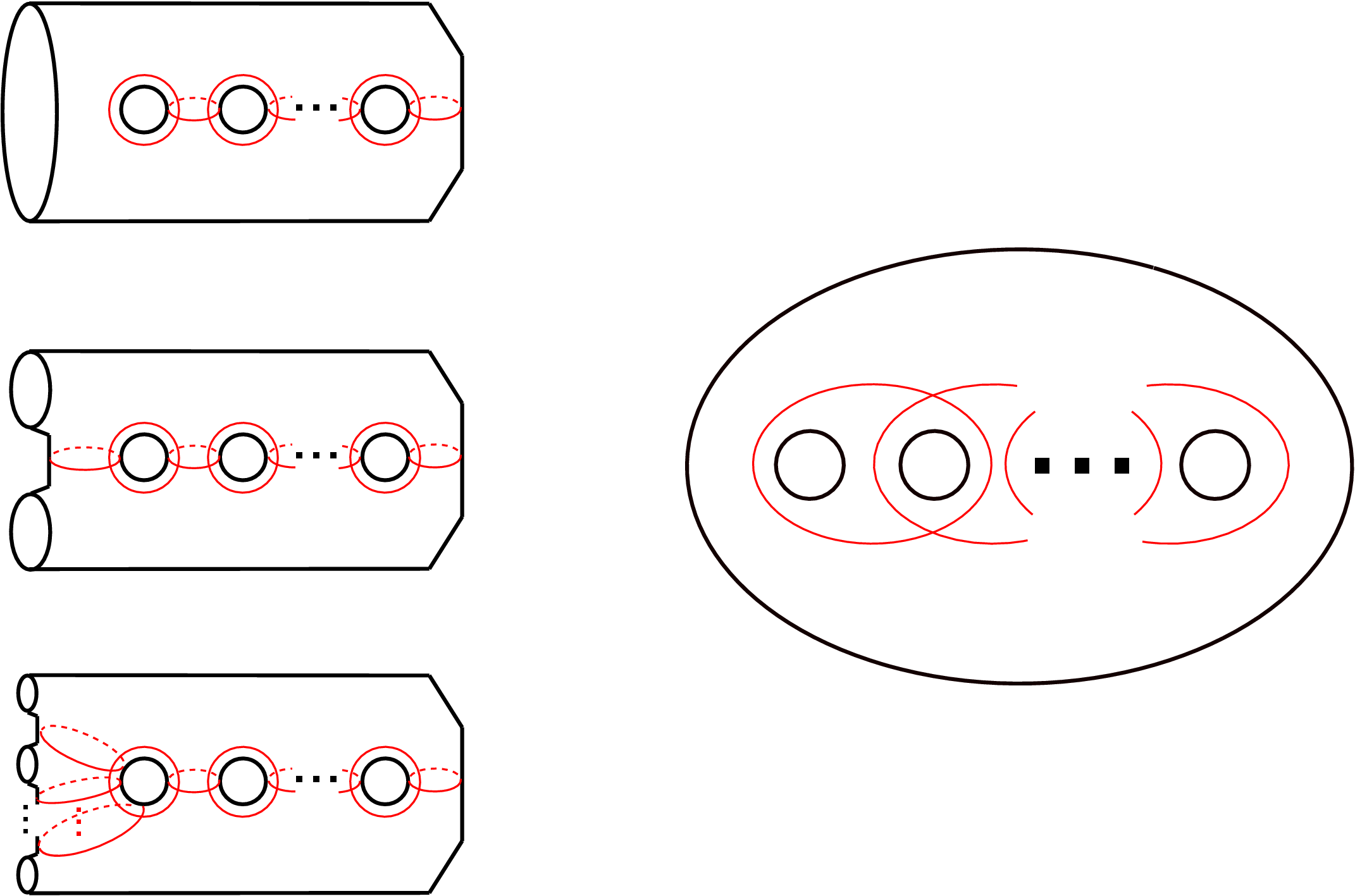}
 \end{center}
 \caption{Curves for $\Sigma_{g,b}$. The left for $g\geq1$ (the top for $b=1$, the middle for $b=2$, the bottom for $b\geq2$.) The right for $g=0$ and $b\geq4$.}
 \label{toannuli}
\end{figure}

We consider the kernel of $[F]$. 
Suppose that $[f]\in Comm([\phi])$ is in the kernel of $[F]$ i.e $F(f)$ is isotopic to the identity map in $S$. 
Note that for every oriented simple cosed curve $l$ in $\Sigma_{g,b}$, $F(l\times \{0\})\subset \Sigma_{g,b}\times \{0\}$ and $l\times \{0\} \subset \Sigma_{g,b}\times \{0\}$ are isotopic in $S$. 
This implies that $F(l\times \{0\})\subset \Sigma_{g,b}\times \{0\}$ and $l\times \{0\} \subset \Sigma_{g,b}\times \{0\}$ are isotopic in $\Sigma_{g,b}\times \{0\}$. 
Thus we can assume that $f$ preserves all oriented simple closed curves in $\Sigma_{g,b}$. 
Take oriented simple closed curves such that they cut $\Sigma_{g,b}$ into annuli whose core curves are some boundary components of $\Sigma_{g,b}$ and disks and a pair of pants each of whose boundary components are some boundary components of $\Sigma_{g,b}$. See Figure~\ref{toannuli}, where orientations are omitted. Note that we need not curves for $(g,b)=(0,1),(0,2),(0,3)$. 
Since $f$ preserves these oriented simple closed curves, $f$ is a composition of Dehn twists along some boundary components of $\Sigma_{g,b}$. 

Conversely, consider $F(t_{c})$ for some boundary component $c$ of $\Sigma_{g,b}$, noting that $t_{c}$ commutes with $\phi$ since $\phi$ fixes the boundary. 
In $\Sigma_{g,b}\times \{1\}\subset S$, $F(t_{c})$ is $t_{c\times \{1\}}$, and in $\Sigma_{g,b}\times \{0\}\subset S$, $F(t_{c})$ is ${t_{c\times \{0\}}}^{-1}$ since the orientation on $\Sigma_{g,b}\times \{0\}$ coming from $S$ is the opposite one to that of coming from $\Sigma_{g,b}$. 
Therefore, as an element of $MCG_{+}(S)$, $[F]([t_{c}])$ is the identity element. 

By the argument above, we can conclude that $Ker([F])=\left< [t_{c}]\mid  c {\rm \ is\ a\ boundary\ component\ of\ }\Sigma_{g,b} \right>$, and this finishes the proof of the former part of Theorem~\ref{thm}.

\section{A proof of the latter part of Theorem~\ref{thm}}\label{latterproof}
In this section, we give a proof of the latter part of Theorem~\ref{thm}. 
Take any orientation reversing involution $\iota$ of $\Sigma_{g,b}$ and an ordered tuple of paiwise disjoint properly embedded oriented arcs in $\Sigma_{g,b}$ such that they cut $\Sigma_{g,b}$ into a disk, $\mathbb{J}=(J_{1},\dots,J_{2g+b-1})$. 
Let ${J_{i}}'$ be an oriented arc on $\Sigma_{g,b}$ obtained by slightly sliding $J_{i}$ to the right side so that the endpoints of ${J_{i}}'$ and $\iota \left({J_{i}}'\right)$ are not on the arcs in $\mathbb{J}$. 
And let ${A_{i}}'$ and ${B_{i}}'$ be oriented simple closed curves on $S$ obtained by replacing all $J_{i}$'s with ${J_{i}}'$'s in the definitions of $A_{i}$ and $B_{i}$, respectively. 
Note that ${A_{i}}'$ and ${B_{i}}'$ are parallel to $A_i$ and $B_i$, respectively. 

Take $[h]\in MCG_{+}(S)$ such that $h$ preserves the binding $\left( \partial \Sigma_{g,b}\times \{\frac{1}{2}\} \right) \subset S$ as a set and reverses the orientation. 
In this situation, we get  two orientation preserving self-homeomorphisms of $\Sigma_{g,b}$ fixing the boundary pointwise $f_{01}$ and $f_{10}$ such that
\begin{itemize}
\item $h\mid_{\partial \Sigma_{g,b}\times [0,1]}(q,t)=\left( \iota(q), 1-t\right)$ for all $q\in \partial \Sigma_{g,b}$ and $t\in[0,1]$. 
\item $h\mid_{\Sigma_{g,b}\times \{0\}}(p,0)=\left( f_{01}\circ \iota(p),1\right)$ for all $p\in \Sigma_{g,b}$.
\item $h\mid_{\Sigma_{g,b}\times \{1\}}(p,1)=\left( f_{10}\circ \iota(p),0\right)$ for all $p\in \Sigma_{g,b}$.
\end{itemize}

We consider conditions for this $h$'s preserving $S_{\mathbb{A}}$ and $S_{\mathbb{B}}$. 

\subsection*{A condition for preserving $S_{\mathbb{A}}$}
For all $i=1,\dots,2g+b-1$,\\
 $h({A_{i}}')=\left( f_{01}\circ \iota({J_{i}}')\times \{1\} \right)\cup \left( \{\iota(t({J_{i}}'))\}\times rev([0,1]) \right)\cup \left(  rev\left( f_{10}\circ \iota ({J_{i}}')\right) \times \{0\} \right)\cup \left( \{\iota(s({J_{i}}'))\}\times [0,1]\right)$. \\
Then $\mathcal{W}_{\mathbb{A}}\left( h({A_{i}}') \right)=Cl\left( {W_{rev(\mathbb{J})}\left( f_{01}\circ \iota ({J_{i}}') \right)} \cdot \overline{W_{\mathbb{J}}\left( rev\left( f_{10}\circ \iota({J_{i}'}) \right)\right)} \right)$. 
Note that the right side of $J_{i}\times \{0\}$ in $\Sigma_{g,b}\times\{0\}$ with the orientation coming from that of $\Sigma_{g,b}$ corresponds to the left side of $A_i$ in $\partial S_{\mathbb{A}}$, and that the right side of $rev(J_{i})\times \{1\}$ in $\Sigma_{g,b}\times\{1\}$ with the orientation coming from that of $\Sigma_{g,b}$, which is the same as that coming from that of $S$, corresponds to the right side of $A_i$ in $\partial S_{\mathbb{A}}$. 
Using Remark~\ref{ref,rev,homeo}, we know that $\overline{W_{\mathbb{J}}\left( rev\left( f_{10}\circ \iota ({J_{i}}')\right) \right)}=\left(W_{rev(\mathbb{J})}\left( f_{10}\circ \iota({J_{i}}') \right) \right)^{-1}$. 
Therefore, $h({A_{i}}')$ bounds a disk in $S_{\mathbb{A}}$ (i.e $\mathcal{W}_{\mathbb{A}}\left( h({A_{i}}') \right)$ is the empty word) if and only if $W_{rev(\mathbb{J})}\left( f_{01}\circ \iota ({J_{i}}') \right)=W_{rev(\mathbb{J})}\left( f_{10}\circ \iota({J_{i}}') \right)$. 
By Remark~\ref{uniqueword} since $f_{01}$ and $f_{10}$ fix the endpoints of $\iota({J_{i}}')$, this is equivalent to that $f_{01}\left(\iota({J_{i}'})\right)$ and $f_{10}\left(\iota({J_{i}}')\right)$ are isotopic in $\Sigma_{g,b}$ relative to their endpoints. 

Since curves in $\mathbb{A}$ cut $S$ into a planar surface, $h$ preserves $S_{\mathbb{A}}$ if and only if $h({A_{i}}')$ bounds a disk in $S_{\mathbb{A}}$ for all $i=1,\dots,2g+b-1$. 
By the argument in previous paragraph, it is equivalent to that $f_{01}\left(\iota({J_{i}'})\right)$ and $f_{10}\left(\iota({J_{i}}')\right)$ are isotopic in $\Sigma_{g,b}$ relative to their endpoints for all $i=1,\dots,2g+b-1$. 
Note that we can isotope $f_{01}\left(\iota({J_{i}'})\right)$ to $f_{10}\left(\iota({J_{i}}')\right)$ for all $i$ simultaneously. 
Since arcs in $\iota(\mathbb{J})$ cut $\Sigma_{g,b}$ into a disk, this is equivalent to that $f_{01}$ and $f_{10}$ are isotopic under the isotopy fixing the boundary.

\subsection*{A condition for preserving $S_{\mathbb{B}}$}
For all $i=1,\dots,2g+b-1$,\\
 $h({B_{i}}')=\left( f_{01}\circ \iota \left( \phi({J_{i}}')\right)\times \{1\} \right)\cup \left( \{\iota \left(t({J_{i}}')\right)\}\times rev([0,1]) \right)\cup \left(  rev\left( f_{10}\circ \iota({J_{i}}')\right) \times \{0\} \right)\cup \left( \{\iota \left(s({J_{i}}')\right) \}\times [0,1] \right)$. \\
Then $\mathcal{W}_{\mathbb{B}}\left( h({B_{i}}') \right)=Cl\left( \overline{W_{rev(\mathbb{J})}\left( f_{01}\circ \iota \circ \phi({J_{i}}') \right) }\cdot W_{\phi(\mathbb{J})}\left( rev\left( f_{10}\circ \iota({J_{i}'}) \right)\right) \right)$. 
Note that the right side of $\phi(J_{i})\times \{0\}$ in $\Sigma_{g,b}\times\{0\}$ with the orientation coming from that of $\Sigma_{g,b}$ corresponds to the right side of $B_i$ in $\partial S_{\mathbb{B}}$, and that the right side of $rev(J_{i})\times \{1\}$ in $\Sigma_{g,b}\times\{1\}$ with the orientation coming from that of $\Sigma_{g,b}$, which is the same as that coming from that of $S$, corresponds to the left side of $B_i$ in $\partial S_{\mathbb{B}}$. 
Using Remark~\ref{ref,rev,homeo}, we know that $\overline{W_{rev(\mathbb{J})}\left( f_{01}\circ \iota \circ \phi({J_{i}}') \right) }=W_{\mathbb{J}}\left( f_{01}\circ \iota \circ \phi({J_{i}}') \right)$ and that $W_{\phi(\mathbb{J})}\left( rev\left( f_{10}\circ \iota({J_{i}'}) \right)\right)=W_{\mathbb{J}}\left( rev\left( \phi^{-1}\circ f_{10}\circ \iota({J_{i}'}) \right)\right)= \left(W_{\mathbb{J}}\left(  \phi^{-1}\circ f_{10}\circ \iota({J_{i}'}) \right)\right)^{-1}$. 
Therefore, $h({B_{i}}')$ bounds a disk in $S_{\mathbb{B}}$ (i.e $\mathcal{W}_{\mathbb{B}}\left( h({B_{i}}') \right)$ is the empty word) if and only if $W_{\mathbb{J}}\left( f_{01}\circ \iota \circ \phi({J_{i}}') \right)=W_{\mathbb{J}}\left(  \phi^{-1}\circ f_{10}\circ \iota({J_{i}'}) \right)$. 
By Remark~\ref{uniqueword} since $f_{01}$, $f_{10}$ and $\phi$ fix the endpoints of $\iota({J_{i}}')$, this is equivalent to that $f_{01}\circ \iota \circ \phi({J_{i}}')$ and $\phi^{-1}\circ f_{10}\circ \iota({J_{i}'}) $ are isotopic in $\Sigma_{g,b}$ relative to their endpoints. 

Since curves in $\mathbb{B}$ cut $S$ into a planar surface, $h$ preserves $S_{\mathbb{B}}$ if and only if $h({B_{i}}')$ bounds a disk in $S_{\mathbb{B}}$ for all $i=1,\dots,2g+b-1$. 
By the argument in previous paragraph, it is equivalent to that $f_{01}\circ \iota \circ \phi({J_{i}}')$ and $\phi^{-1}\circ f_{10}\circ \iota({J_{i}'}) $ are isotopic in $\Sigma_{g,b}$ relative to their endpoints for all $i=1,\dots,2g+b-1$. 
Note that we can isotope $f_{01}\circ \iota \circ \phi({J_{i}}')$ to $\phi^{-1}\circ f_{10}\circ \iota({J_{i}'}) $ simultaneously. 
Since arcs in $\mathbb{J}$ cut $\Sigma_{g,b}$ into a disk, this is equivalent to that $f_{01}\circ \iota \circ \phi$ and $\phi^{-1}\circ f_{10}\circ \iota$ are isotopic under the isotopy fixing the boundary.  
Composing the involution $\iota$ from the right, we get the equation $\phi^{-1}\circ f_{10}=f_{01}\circ \iota \circ \phi \circ \iota$. 

\vspace{1.0cm}

By the arguments above, we conclude that $\mathcal{G}$ has an element one of whose representatives fixes the binding as a set and reverses the orientation of the binding if and only if there exists $[f]\in MCG_{+}(\Sigma_{g,b})$ such that $\phi^{-1}\circ f=f\circ \iota \circ \phi \circ \iota$. This finishes the proof of the latter part of Theorem~\ref{thm}.

\section{An example} \label{example}
In this section, we give some example of a Goeritz group of a Heegaard splitting induced by an openbook decomposition. 

\begin{figure}[htbp]
 \begin{center}
  \includegraphics[width=40mm]{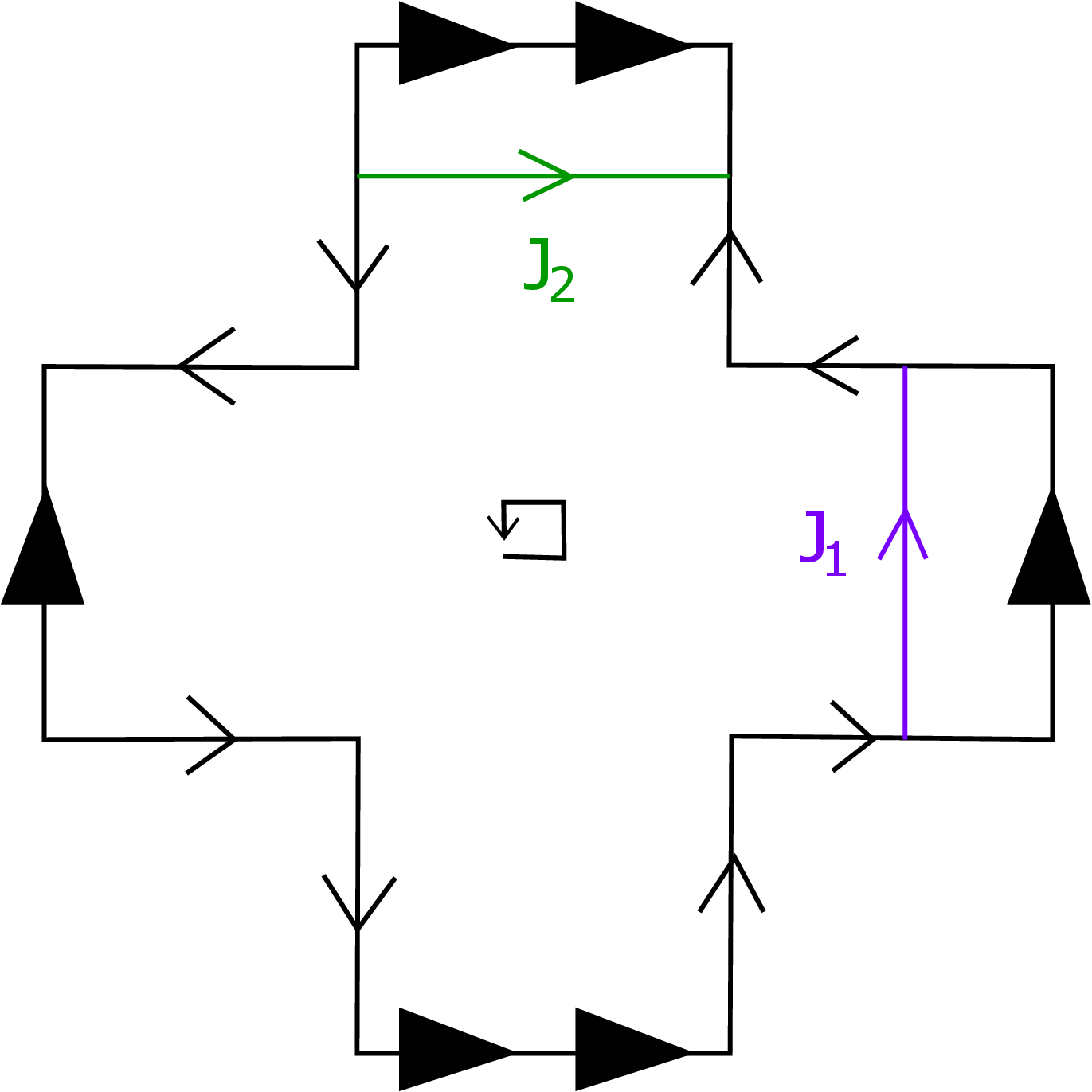}
 \end{center}
 \caption{Arcs $J_1$ and $J_2$ cutting $\Sigma_{1,1}$ into a disk.}
 \label{arctorus}
\end{figure}

\begin{figure}[htbp]
 \begin{center}
  \includegraphics[width=100mm]{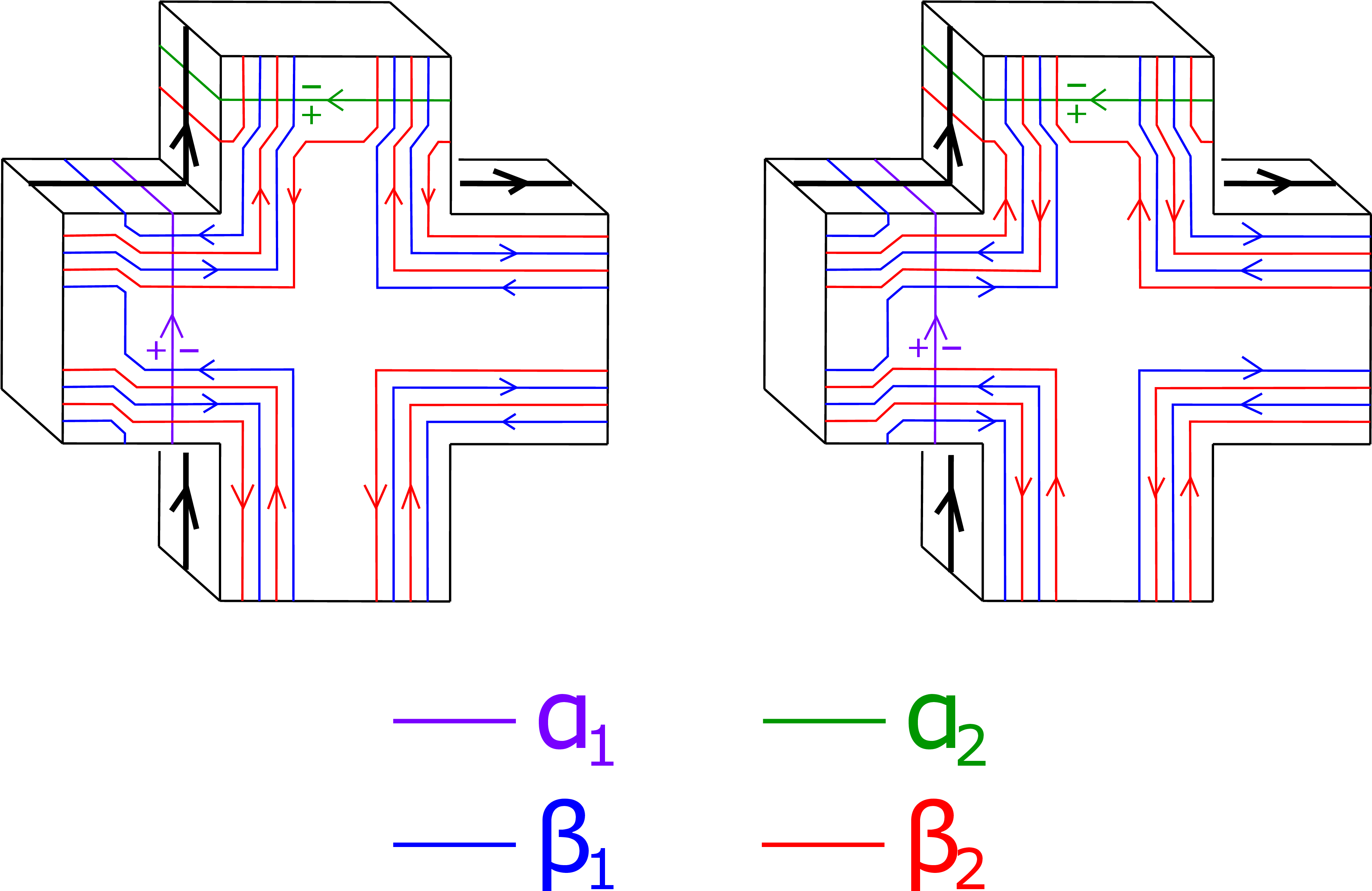}
 \end{center}
 \caption{Diagrams for $n=1$ (left), for $n=-1$ (right).}
 \label{diagram_example}
\end{figure}

In the following, let $N$ be a 3-manifold which admits an openbook decomposition whose page is homeomorphic to $\Sigma_{1,1}$ and whose monodromy is ${t_{\partial}}^n$ for a non-zero integer $n$, where $t_{\partial}$ is the right-hand Dehn twist along curve obtained by slightly pushing $\partial \Sigma_{1,1}$ into the interior of $\Sigma_{1,1}$. 
By fixing properly embedded arc $J_{1}$ and $J_{2}$ as in Figure~\ref{arctorus}, where the arrow in the center indicates the orientation of $\Sigma_{1,1}$, we can construct a Heegaard diagram $\left( S;(\alpha_{1},\alpha_{2}), (\beta_{1},\beta_{2}) \right)$ as in Subsection~\ref{hdiagram}. 
We give figures of $\left( S;(\alpha_{1},\alpha_{2}), (\beta_{1},\beta_{2}) \right)$ for $n=\pm1$ in Figure~\ref{diagram_example}, where to obtain $S$, we reverse $\Sigma_{1,1}$ and thicken it to the back side of the paper. 
Thick black curve presents the binding. 
 
Let $\mathcal{G}$ be the Goeritz group of this Heegaard splitting, and $\mathcal{G}_{bind}$ a subgroup of $\mathcal{G}$ consisting of elements which preserve the binding and its orientation.
We will prove the following: 
\begin{prop} \label{prop}
$\mathcal{G}=\mathbb{Z}/2\mathbb{Z}\ltimes \mathcal{G}_{bind}\cong \left<A,B,r \mid ABA=BAB, (AB)^{6}=1 , r^{2}=1, rAr=A^{-1}, rBr=B^{-1}\right>$.
\end{prop}

We prove Proposition~\ref{prop} in three steps:
\begin{itemize}
\item We see that every element of $\mathcal{G}$ fixes the binding as a set (Subsection~\ref{bindingfix}). 
\item We give an element $r$ of order $2$ of $\mathcal{G}$ reversing the binding (Subsection~\ref{bindingreverse}).
\item We directly compute $\mathcal{G}_{bind}$, and the conjugate action of $r$ to $\mathcal{G}_{bind}$ (Subsection~\ref{presentation}). 
\end{itemize}

\begin{rmk}
For the case where $n=0$, namely the monodromy is the identity map, we get a Heegaard splitting of genus two of $\left(S^{2}\times S^{1}\right) \# \left(S^{2}\times S^{1}\right)$. 
It is known that the Goeritz group of this splitting is isomorphic to the handlebody group of genus two. 
\end{rmk}

\subsection{Simple closed curves on a closed surface of genus two and GOF-knots on the boundary of a handlebody of genus two}
In this subsection, we recall some definitions about simple closed curves on a closed surface of genus two and GOF-knots on the boundary of a handlebody of genus two, which is appeared in \cite{ck2}, \cite{ck3}, \cite{sekino1}, and will be used for the proof of Proposition~\ref{prop}. 
All curves in a surface is assumed that they intersect each other minimally and transversely. 

\subsubsection{Simple closed curves on a closed surface of genus two}
Let $\Sigma$ be a connected oriented closed surface of genus two, and $(l_{1},l_{2})$ an ordered pair of pairwise disjoint oriented simple closed curves in $\Sigma$ cutting $\Sigma$ into a planar surface. 
Let $P\left(\Sigma;l_1,l_2\right)$ denote the closure of the result of cutting $\Sigma$ along $l_1$ and $l_2$. 
We give $P\left(\Sigma;l_1,l_2\right)$ the orientation coming from $\Sigma$. 
For each $i=1,2$, ${l_{i}}^{+}$ denotes the boundary component coming from $l_{i}$ and whose orientation as a subset of $P\left(\Sigma;l_1,l_2\right)$ is the same as that of coming from $l_{i}$, and ${l_{i}}^{-}$ denotes the other boundary component coming from $l_{i}$. 
For every multi-curve in $\Sigma$, this is presented as a set of properly embedded arcs or simple closed curves in $P\left(\Sigma;l_1,l_2\right)$, we call it a {\it planer diagram} of the multi-curve in $P\left(\Sigma;l_1,l_2\right)$. 

In the following definitions and facts, let $l$ be an oriented simple closed curve in $\Sigma$.
\begin{defini}
\begin{itemize}
\item We assign a {\it cyclic sequence of letters} $s_{(l_1,l_2)}(l)$ in alphabets $\{x_{l_1},x_{l_2}\}$ as follows: 
Follow $l$ from any point, and we start with the empty sequence. If we hit $l_{i}$ from the right side (or left side), then we add $x_{l_i}$ (or ${x_{l_i}}^{-1}$, respectively) to the right of the sequence we have. 
When we return to the starting point, connect the last of the sequence to the initial of it to get a cyclic sequence. 
Note that a sequence of letters is not a word, and we do not reduce a part such as $x_{l_i}{x_{l_{i}}}^{-1}$. 
\item We assign a {\it cyclic word} $w_{(l_1,l_2)}(l)$ in alphabets $\{x_{l_1},x_{l_2}\}$ by reducing $s_{(l_1,l_2)}(l)$ as cyclic words. 
\end{itemize}
\end{defini}

\begin{defini}
Suppose $s_{(l_1,l_2)}(l)$ has one of $x_{l_1}x^{-1}_{l_1}$, $x^{-1}_{l_1}x_{l_1}$, $x_{l_2}x^{-1}_{l_2}$ and $x^{-1}_{l_2}x_{l_2}$, say $x_{l_1}x^{-1}_{l_1}$. 
This subsequence of letters corresponds to a subarc $c$ of $l$ starting from and ending on $l^{+}_{1}$ in the planar diagram $P\left(\Sigma;l_1,l_2\right)$. 
This $c$ separates other three boundary components of $P\left(\Sigma;l_1,l_2\right)$ into two sets, none of which are empty since $l$ intersects with $l_1$ essentially. 
Thus there is an oriented arc $\alpha$ on $P\left(\Sigma;l_1,l_2\right)$ starting from $l^{+}_{1}$ and ending on another boundary component, denoted by $A$, such that $c$ is the boundary of a small neighborhood of $ \alpha \cup A$ in the interior of $P\left(\Sigma;l_1,l_2\right)$. 
This $\alpha$ determines $c$ except for the orientation, and is called the {\it corresponding arc} of $c$. 
This $A$ is not $l^{-}_{1}$ since otherwise $|l \cap l^{+}_{1}|$ were grater than $|l \cap l^{-}_{1}|$. 
Note that if there is another subarc $c'$ of $l$ representing $x_{l_1}x^{-1}_{l_1}$, then its corresponding arc $\alpha'$ is isotopic (endpoints can move on cut ends) to $\alpha$. 
This is because $c$ and $c'$ would intersect otherwise.
\end{defini}

\begin{fact}\label{samenumber}(\cite{sekino1})\\
For each $i=1,2$, the number of arcs in a diagram of $l$ in $P\left(\Sigma;l_1,l_2\right)$ starting from and ending on ${l_i}^{+}$ is the same as that of arcs in a diagram of $l$ in $P\left(\Sigma;l_1,l_2\right)$ starting from and ending on ${l_i}^{-}$. 
\end{fact}

\begin{fact}\label{condition_reduced}(Lemma 3.3. of \cite{ck3})\\
If $s_{(l_1,l_2)}(l)$ has one of $x_{l_1}{x_{l_2}}^{m}{x_{l_1}}^{-1}$, ${x_{l_1}}^{-1}{x_{l_2}}^{m}x_{l_1}$, $x_{l_2}{x_{l_1}}^{m}{x_{l_2}}^{-1}$ and ${x_{l_2}}^{-1}{x_{l_1}}^{m}x_{l_2}$ for a non-zero integer $m$, then $s_{(l_1,l_2)}(l)=w_{(l_1,l_2)}(l)$ i.e. $s_{(l_1,l_2)}(l)$ is cyclically reduced as a word. 
\end{fact}

\subsubsection{GOF-knot of a handlebody of genus two}
Let $H$ be a handlebody of genus two, oriented as a subset of $\mathbb{R}^3$. 
In \cite{sekino1}, an oriented simple closed curve $l$ in $\partial H$ is called a {\it GOF-knot} of $H$ if there is an homeomorphism between $\left( \Sigma_{1,1}\times[0,1], \partial \Sigma_{1,1}\times \{\frac{1}{2}\} \right)$ and $(H,l)$ as pair of unoriented manifolds. 
Note that if $l$ is a GOF-knot of $H$, then $l$ with the opposite orientation is also a GOF-knot of $H$. 

Take an ordered pair of oriented properly embedded pairwise disjoint disks $(D_{1},D_{2})$ such that $D_1$ and $D_2$ cut $H$ into a ball. 
Then we have the following. 
\begin{fact}\label{gof}
An oriented simple closed curve $l$ in $\partial H$ is a GOF-knot of $H$ if and only if the cyclic word $w_{(\partial D_1,\partial D_2)}$ is $[[x_{\partial D_1}, x_{\partial D_2}]]$ or $[[x_{\partial D_1}, {x_{\partial D_2}}^{-1}]]$. 
\end{fact}

Note that in our example $\left( S; (\alpha_1,\alpha_2),(\beta_1,\beta_2)\right)$, the binding is a GOF-knot of both handlebodies of genus two $S_{(\alpha_1,\alpha_2)}$ and $S_{(\beta_1,\beta_2)}$.

\subsection{The binding is fixed as a set} \label{bindingfix}
Take an element $[h]\in \mathcal{G}$, regarded as the mapping class of $S$. 
Let $l$ denote the binding. 
Since $l$ is a GOF-knot of both handlebodies of genus two and $h$ fixes the Heegaard splitting, $h(l)$ is also a GOF-knot of both handlebodies of genus two. 
By perturbing $h$, we assume that $h(l)$ intersects $\alpha_1$, $\alpha_2$, $\beta_1$ and $\beta_2$ transversely and minimally. 

\begin{figure}[htbp]
 \begin{center}
  \includegraphics[width=80mm]{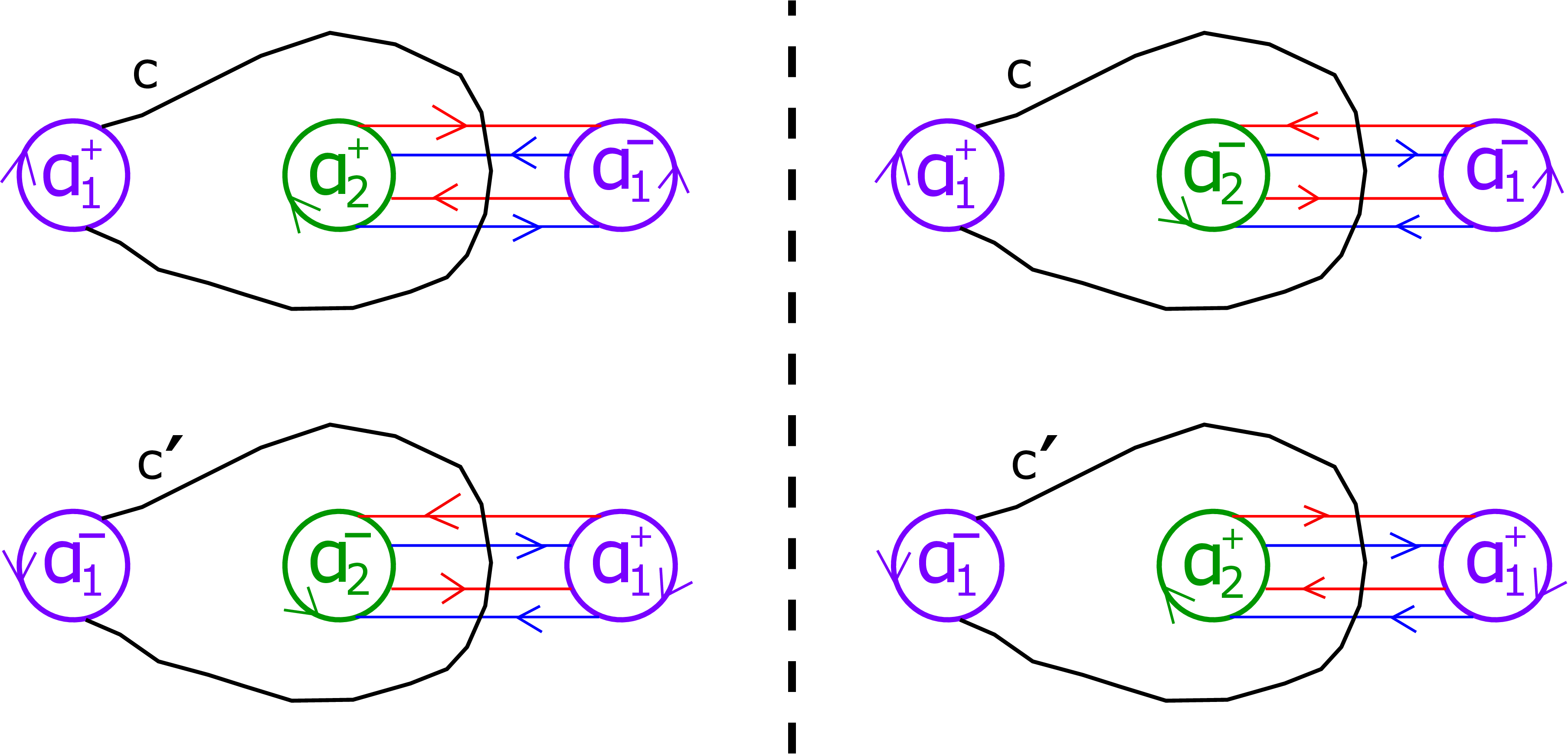}
 \end{center}
 \caption{The intersections of $c$ and $c'$ for $n=1$, the case is divided by the terminal point of the corresponding arc of $c$.}
 \label{to_reduced}
\end{figure}

If a diagram of $h(l)$ in $P\left(S;\alpha_1,\alpha_2\right)$ has an arc $c$ starting from and ending on ${\alpha_1}^{+}$, then the diagram also has an arc $c'$ starting from and ending on ${\alpha_1}^{-}$ by Fact~\ref{samenumber}. 
Conversely, if a diagram of $h(l)$ in $P\left(S;\alpha_1,\alpha_2\right)$ has an arc $c'$ starting from and ending on ${\alpha_1}^{-}$, then the diagram also has an arc $c$ starting from and ending on ${\alpha_1}^{+}$ by Fact~\ref{samenumber}. 
In this case, we see that $s_{(\beta_1,\beta_2)}\left(h(l)\right)=w_{(\beta_1,\beta_2)}\left(h(l)\right)$ by Fact~\ref{condition_reduced} and $h(l)$ hits $\beta_{1}\cup \beta_2$ at least eight times. See Figure~\ref{to_reduced}. 
This contradicts to that $h(l)$ is a GOF-knot of $S_{(\beta_1,\beta_2)}$ by Fact~\ref{gof}. 
Similarly, a diagram of $h(l)$ in $P\left(S;\alpha_1,\alpha_2\right)$ has neither arcs starting from and ending on ${\alpha_2}^{+}$ nor arcs starting from and ending on ${\alpha_2}^{-}$. 
Thus we know that $s_{(\alpha_1,\alpha_2)}\left(h(l)\right)=w_{(\alpha_1,\alpha_2)}\left(h(l)\right)$.

By Fact~\ref{gof}, the diagram of $h(l)$ in $P\left( S;\alpha_1,\alpha_2\right)$ consists of four arcs. 
As unoriented arcs, they are one connecting ${\alpha_{1}}^{+}$ and ${\alpha_2}^{-}$, one connecting ${\alpha_2}^{+}$ and ${\alpha_1}^{+}$, one connecting ${\alpha_1}^{-}$ and ${\alpha_2}^{+}$ and one connecting ${\alpha_2}^{-}$ and ${\alpha_1}^{-}$. 
If the arc of the diagram of $h(l)$ in $P\left( S;\alpha_1,\alpha_2\right)$ connecting ${\alpha_{1}}^{+}$ and ${\alpha_2}^{-}$ is not parallel to the arcs of the diagram of $\beta_{1}\cup \beta_2$ in $P\left( S;\alpha_1,\alpha_2\right)$ connecting ${\alpha_{1}}^{+}$ and ${\alpha_2}^{-}$, then this arc passes over the arcs of the diagram of $\beta_{1}\cup \beta_2$ in $P\left( S;\alpha_1,\alpha_2\right)$ connecting ${\alpha_{1}}^{-}$ and ${\alpha_2}^{+}$. 
In this case, the arc of the diagram of $h(l)$ in $P\left( S;\alpha_1,\alpha_2\right)$ connecting ${\alpha_{1}}^{-}$ and ${\alpha_2}^{+}$ passes over the arcs of the diagram of $\beta_{1}\cup \beta_2$ in $P\left( S;\alpha_1,\alpha_2\right)$ connecting ${\alpha_{1}}^{+}$ and ${\alpha_2}^{-}$. 
Then by Fact~\ref{condition_reduced}, $s_{(\beta_1,\beta_2)}\left(h(l)\right)=w_{(\beta_1,\beta_2)}\left(h(l)\right)$ and $h(l)$ hits $\beta_{1}\cup \beta_2$ at least eight times. 
This contradicts to that $h(l)$ is a GOF-knot of $S_{(\beta_1,\beta_2)}$ by Fact~\ref{gof}. 
Thus the arc of the diagram of $h(l)$ in $P\left( S;\alpha_1,\alpha_2\right)$ connecting ${\alpha_{1}}^{+}$ and ${\alpha_2}^{-}$ is parallel to the arcs of the diagram of $\beta_{1}\cup \beta_2$ in $P\left( S;\alpha_1,\alpha_2\right)$ connecting ${\alpha_{1}}^{+}$ and ${\alpha_2}^{-}$. 
Similarly, the other three arcs of $h(l)$ are parallel to some of arcs in the diagram of $\beta_{1}\cup \beta_2$ in $P\left( S;\alpha_1,\alpha_2\right)$. 

Give $h(l)$ temporary the orientation so that the arc of the diagram of $h(l)$ in $P\left( S;\alpha_1,\alpha_2\right)$ connecting ${\alpha_{1}}^{+}$ and ${\alpha_2}^{-}$ starts from ${\alpha_{1}}^{+}$ and ends on ${\alpha_2}^{-}$. 
Then the (cyclic) order of the connection of four arcs to construct $h(l)$ is the one from ${\alpha_1}^{+}$ to ${\alpha_2}^{-}$ first, the one from ${\alpha_2}^{+}$ to ${\alpha_1}^{+}$ second, ${\alpha_1}^{-}$ to ${\alpha_2}^{+}$ third, ${\alpha_2}^{-}$ to ${\alpha_1}^{-}$ fourth. 
This order is the same as that of $l$. 
Thus, $h(l)$ or $h(l)$ with the opposite orientation is isotopic to ${t_{\alpha_1}}^{u}\circ {t_{\alpha_2}}^{v}(l)$ for some integers $u,v$. 
Note that $w_{(\beta_1,\beta_2)}\left( {t_{\alpha_1}}^{u}\circ {t_{\alpha_2}}^{v}(l)\right)=
\left[ \left[ \left( \left([[x_{\beta_2}, x_{\beta_1}]]\right)^{n} \left([[{x_{\beta_2}}^{-1}, x_{\beta_1}]] \right)^{n} \right)^{u} x_{\beta_1} , \left( \left([[{x_{\beta_1}}^{-1}, x_{\beta_2}]]\right)^{n} \left([[x_{\beta_1}, x_{\beta_2}]] \right)^{n} \right)^{v} x_{\beta_2} \right] \right] $. 
And this cannot be a GOF-knot of $S_{(\beta_1,\beta_2)}$ by Fact~\ref{gof} unless $u=v=0$. 
Hence we conclude that $h(l)$ or $h(l)$ with the opposite orientation is isotopic to $l$ in $S$.

\subsection{An element of $\mathcal{G}$ reversing the binding}\label{bindingreverse}
\begin{figure}[htbp]
 \begin{center}
  \includegraphics[width=60mm]{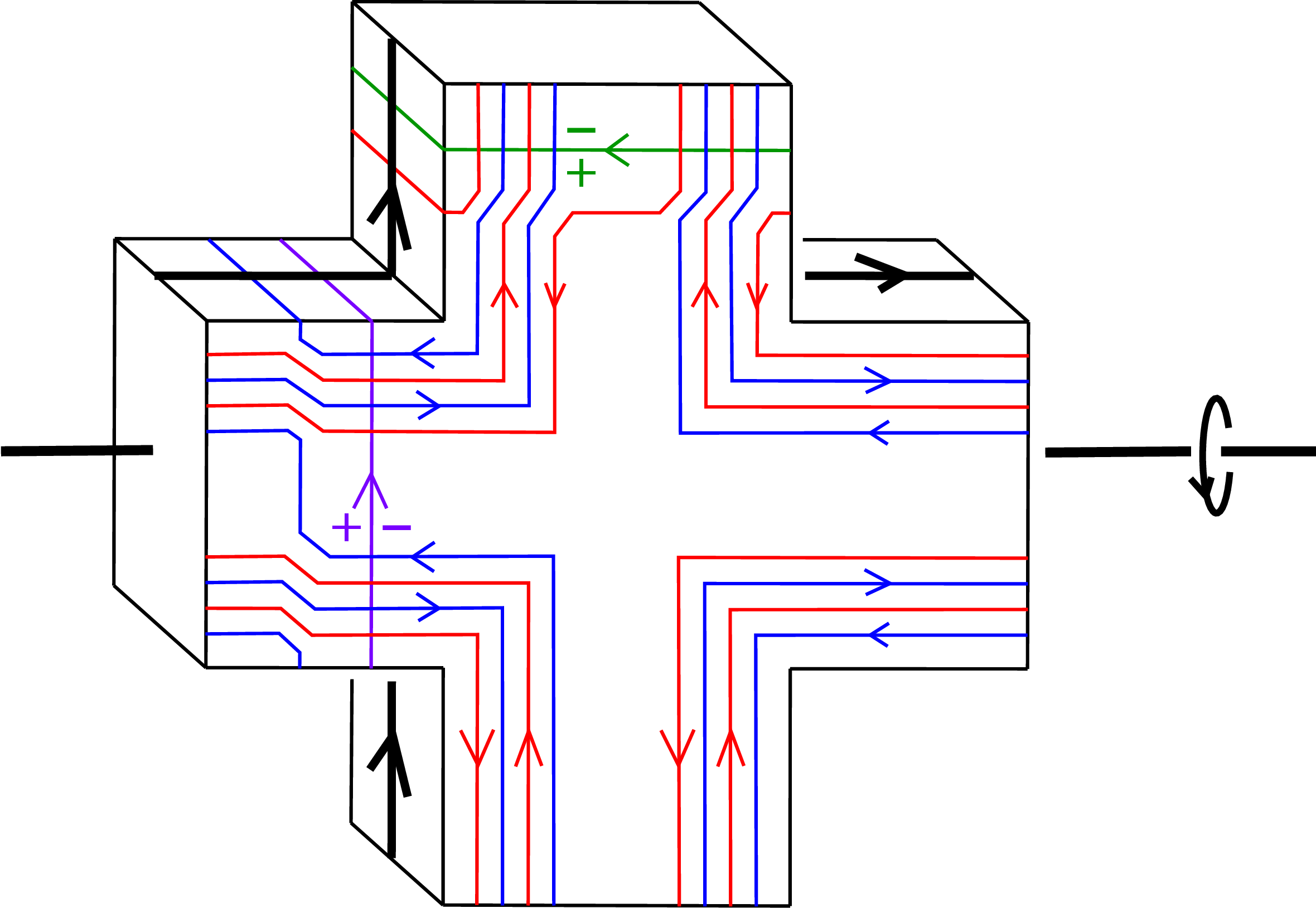}
 \end{center}
 \caption{The homeomorphism $r$ for $n=1$.}
 \label{reversemap}
\end{figure}

Let $r$ be a orientation preserving homeomorphism of $S$ obtained by $\pi$-rotation along the horizontal axis in Figure~\ref{reversemap}. 
Since $r(\alpha_1)=\alpha_1$ and $r(\alpha_2)$ is isotopic to $\alpha_2$ with the opposite orientation, $r$ extends to a self-homeomorphism of $S_{(\alpha_1,\alpha_2)}$. 
And since $r(\beta_1)$ is isotopic to $\beta_1$ and $r(\beta_2)$ is isotopic to $\beta_2$ with the opposite orientation, $r$ extends to a self-homeomorphism of $S_{(\beta_1,\beta_2)}$. 
Thus we know that $r\in \mathcal{G}$. 
Note that $r$ reverses the binding, and $r^{2}$ is the identity map. 

\subsection{A presentation of $\mathcal{G}$}\label{presentation}
Since every element of $\mathcal{G}$ fixes the binding as a set and the element $[r]$ reversing the binding, $\mathcal{G}$ is generated by the element of $\mathcal{G}_{bind}$ and $[r]$. 
Moreover, since $\mathcal{G}_{bind}$ is a normal subgroup of $\mathcal{G}$ and $[r]^{2}=1_{\mathcal{G}}$, we see that $\mathcal{G}=\mathbb{Z}/2\mathbb{Z}\ltimes \mathcal{G}_{bind}$. 
We give a presentation of $\mathcal{G}_{bind}$ and consider the conjugate action of $[r]$. 

\subsubsection{A presentation of $\mathcal{G}_{bind}$}
\begin{figure}[htbp]
 \begin{center}
  \includegraphics[width=40mm]{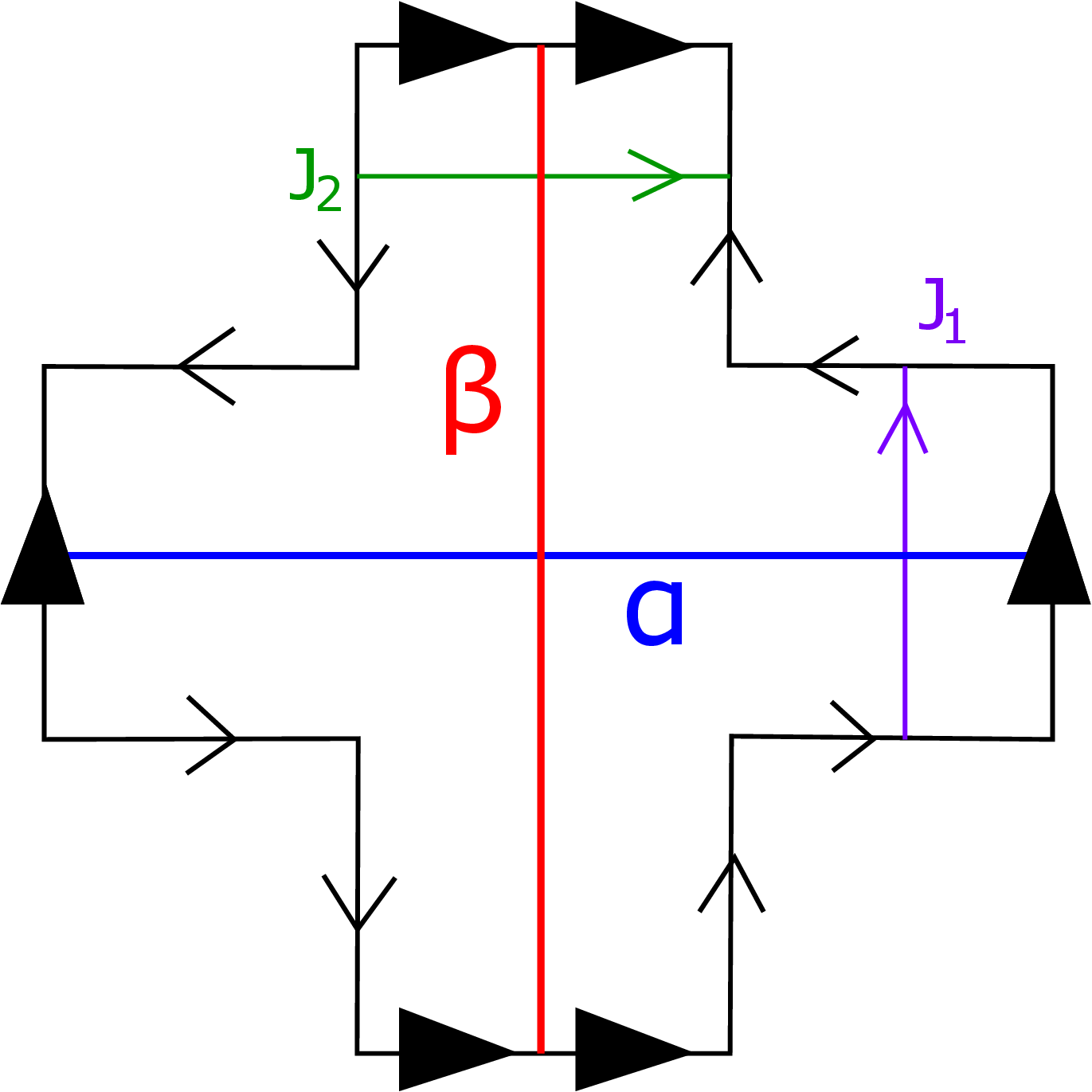}
 \end{center}
 \caption{Two simple closed curves $\alpha$ and $\beta$.}
 \label{alphabeta}
\end{figure}

Since every element of $MCG_{+}(\Sigma_{1,1})$ commutes with the monodromy ${t_{\partial \Sigma_{1,1}}}^{n}$, we know that $\mathcal{G}_{bind}\cong  MCG_{+}(\Sigma_{1,1}) /\left< t_{\partial \Sigma_{1,1}}\right> $ by Theorem~\ref{thm}. 
Fix two unoriented simple closed curves $\alpha$ and $\beta$ in $\Sigma_{1,1}$ as in Figure~\ref{alphabeta}. 
It is known that $MCG_{+}(\Sigma_{1,1})$ has a presentation $\left< A, B\mid ABA=BAB  \right>$, where $A$ and $B$ correspond to $[t_{\alpha}]$ and $[t_{\beta}]$, respectively, see the next paragraph of Theorem 3.14 of  \cite{fm}. 
Under this presentation, $[t_{\partial \Sigma_{1,1}}]$ corresponds to $(AB)^{6}$. 
Hence $\mathcal{G}_{bind}$ has a presentation $\left< A, B \mid ABA=BAB, (AB)^{6}=1 \right>$. 

\subsubsection{A conjugate action of $[r]$}
Recall that $S$ is $\left( \Sigma_{1,1}\times \{0\} \right) \cup \left( \partial \Sigma_{1,1} \times [0,1]\right) \cup \left( \Sigma_{1,1} \times \{1\} \right)$, and that $F\left( t_{\gamma} \right)$ is a composition of ${t_{\gamma \times \{0\}}}^{-1}$ and $t_{\gamma \times \{1\}}$ for a simple closed curve $\gamma$ in $\Sigma_{1,1}$. 
Then we see that $r\circ F(t_{\alpha}) \circ r=F({t_{\alpha}}^{-1})$ and $r\circ F(t_{\beta}) \circ r=F({t_{\beta}}^{-1})$. 
Hence we conclude that $\mathcal{G}\cong \left< A, B, r \mid ABA=BAB, (AB)^{6}=1, r^{2}=1, rAr=A^{-1}, rBr=B^{-1}\right>$.

\vspace{0.5cm}

\ GRADUATE SCHOOL OF MATHEMATICAL SCIENCES, THE UNIVERSITY OF TOKYO, 3-8-1 KOMABA, MEGURO--KU, TOKYO, 153-8914, JAPAN\\
\ \ E-mail address: \texttt{sekinonozomu@g.ecc.u-tokyo.ac.jp}


\begin{thebibliography}{50}
\bibitem[1]{akbas} E. Akbas. A presentation for the automorphisms of the 3-sphere that preserve a genus two Heegaard splitting. Pacific J. Math. 236 (2008), no. 2, 201--222.

\bibitem[2]{alexander} J. W. Alexander. A lemma on systems of knotted curves. Proceedings of the National Academy of Sciences, 9(3):93--95, 1923.

\bibitem[3]{cho1} S. Cho. Homeomorphisms of the 3-sphere that preserve a Heegaard splitting of genus two. Proc. Amer. Math. Soc. 136 (2008), no. 3, 1113--1123.

\bibitem[4]{cho2} S. Cho. Genus two Goeritz groups of lens spaces. Pacific J.Math. 265 (2013), no. 1, 1--16.

\bibitem[5]{ck1} S. Cho, Y. Koda. The genus two Goeritz group of $S^2\times S^1$. Math. Res. Lett. 21 (2014), no. 3, 449--460.

\bibitem[6]{ck2} S. Cho, Y. Koda. Disk complexes and genus two Heegaard splittings for non-prime 3-manifolds. Int. Math. Res. Not. IMRN 2015 (2015), 4344--4371.

\bibitem[7]{ck3} S. Cho, Y. Koda. Connected primitive disk complexes and genus two Goeritz groups of lens spaces. Int. Math. Res. Not. IMRN 2016 (2016), 7302--7340.

\bibitem[8]{ck4} S. Cho, Y. Koda. The mapping class groups of reducible Heegaard splittings of genus two. Transactions of American Mathematical Society 371 (2019), no. 4, 2473--2502.

\bibitem[9]{fm} B. Farb, D. Margalit. A primer on mapping class groups, Princeton Mathematical Series, 49. Princeton University Press, Princeton, NJ, 2012.

\bibitem[10]{fs} M. Freedman, M. Scharlemann. Powell moves and the Goeritz group. arXiv:1804.05909.

\bibitem[11]{goeritz} L. Goeritz. Die Abbildungen der Berzelfl\''{a}che und der Volbrezel vom Gesschlect 2. Abh. Math. Sem. Univ. Hamburg 9 (1933), 244--259.

\bibitem[12]{hempel} J. Hempel. 3-manifolds as viewed from the curve complex. Topology 40 (2001), no. 3, 631--657.

\bibitem[13]{iguchi} D. Iguchi. Thick isotopy property and the mapping class groups of Heegaard splittings. arXiv:2008.11548.

\bibitem[14]{ik1} D. Iguchi, Y. Koda. Twisted book decompositions and the Goeritz groups. Topology Appl. 272 (2020), 107064, 15 pp.

\bibitem[15]{jaco} W. Jaco. Lectures on three-manifold topology. CBMS Regional Conference Series in Mathematics, vol. 43, American Mathematical Society, Providence, RI, 1980.

\bibitem[16]{johnson1} J. Johnson. Mapping class groups of medium distance Heegaard splittings. Proc. Amer. Math. Soc. 138 (2010), no. 12, 4529--4535.

\bibitem[17]{johnson2} J. Johnson. Heegaard splittings and open books. arXiv:1110.2142.

\bibitem[18]{ms} Y. Moriah, J. Schultens. Irreducible Heegaard splittings of Seifert fibered spaces are either vertical or horizontal. Topology 37 (1998), no. 5, 1089--1112.

\bibitem[19]{scharlemann} M. Scharlemann. Automorphisms of the 3-sphere that preserve a genus two Heegaard splitting. Bol. Soc. Mat. Mexicana (3) 10 (2004), Special Issue, 503--514.

\bibitem[20]{sekino1} N. Sekino. Genus one fibered knots in 3-manifolds with reducible genus two Heegaard splittings. Topology Appl. 239 (2018), 46--64.

\bibitem[21]{sekino} N. Sekino. The Goeritz group of a Heegaard splitting of genus two of a Seifert manifold whose base orbifold is sphere with three exceptional points of sufficiently complex coefficients. arXiv:2202.05064.

%\bibitem[1]{quadra} A. Baker. A concise introduction to the theory of numbers. Cambridge University Press, Cambridge, 1984.

%\bibitem[2]{homcob} S. Garoufalidis and J. Levine. Tree-level invariants of three-manifolds, Massey products and the Johnson homomorphism. In Graphs and patterns in mathematics and theoretical physics, volume 73 of Proc. Sympos. Pure Math., pages 173--203. Amer. Math. Soc., Providence, RI, 2005.

%\bibitem[3]{homfib} H. Goda and T. Sakasai. Homology cylinders and sutured manifolds for homologically fibered knots. Tokyo J. Math., 36(1):85--111, 2013.

%\bibitem[4]{density} S. Lang. Algebraic number theory, volume 110 of Graduate Texts in Mathematics. Springer-Verlag, New York, second edition, 1994. 

%\bibitem[5]{lickorish} W. B. R. Lickorish, An introduction to knot theory, Graduate Texts in Mathematics, vol. 175, Springer-Verlag, New York, 1997.

%\bibitem[6]{conti} E. Manfredi, A. Savini. Fibered knots and links in lens spaces. arXiv:1502.03345.

%\bibitem[7]{nozaki} Y. Nozaki. Every lens space contains a genus one homologically fibered knot. Illinois J. Math. volume 62, Number 1--4 (2018), 99--111.

%\bibitem[8]{hc} T. Sakasai. Johnson-Morita theory in mapping class groups and monoids of homology cobordisms of surfaces. Winter Braids Lecture Notes, 3(4):1--25, 2016.

%\bibitem[9]{cycle} S. H. Weintraub. Galois theory. Universitext. Springer, New York, second edition, 2009.

\end{thebibliography}
\end{document}